\documentclass[11pt]{article}

\renewcommand\appendix{
  \setcounter{section}{0}%
  \setcounter{subsection}{0}%
\renewcommand{\theequation}{A.\arabic{equation}} 
 \gdef\thesection{A}\section
}
\renewcommand{\theequation}{\arabic{section}.\arabic{equation}}
\usepackage{amsmath}
\usepackage{amssymb}
\usepackage{amscd}
\newcommand{\KK}{{\rm K}}
\newcommand{\Pic}{{\rm Pic}}
\newcommand{\FR}{{\rm FR}}
\newcommand{\I}{{\rm I}}

\newcommand{\Ball}{{\rm B}}

\renewcommand{\P}{\mathbb{P}}
\newcommand{\C}{\mathbb{C}}
\newcommand{\B}{{\cal B}}

\newcommand{\Z}{\mathbb{Z}}
\newcommand{\id}{{\rm id}}

\renewcommand{\H}{{\cal H}}
\newcommand{\G}{{\cal G}}

\newcommand{\aire}{{\rm aire}}
\newcommand{\trace}{{\rm trace\ \!}}

\newcommand{\ddc}{{\rm dd^c}}
\renewcommand{\d}{{\rm d}}

\newcommand{\dist}{{\rm dist}}
\newcommand{\diam}{{\rm diam\ \!}}
\newcommand{\FS}{{\rm FS}}

\newcommand{\F}{{\cal F}}

\newcommand{\E}{{\cal E}}

\newcommand{\mult}{{\rm mult}}

\newcommand{\voir}{{\it voir }}
\renewcommand{\Re}{{\rm Re}}

\renewcommand{\O}{{\rm O}}
\renewcommand{\o}{{\rm o}}

\newtheorem{theoreme}{Th\'eor\`eme}[section]
\newtheorem{proposition}[theoreme]{Proposition}
\newtheorem{corollaire}[theoreme]{Corollaire}
\newtheorem{lemme}[theoreme]{Lemme}

\newtheorem{exemple}[theoreme]{Exemple}

\newtheorem{remarque}[theoreme]{Remarque}

\newenvironment{preuve}{\begin{trivlist}
\item[]{\bf D\'emonstration.}}
{\par\hfill $\square$\end{trivlist}}
\title{Suites d'applications m\'eromorphes multivalu\'ees\\ 
et courants laminaires} 
\author{Tien-Cuong Dinh}
\begin{document}
\maketitle
\begin{abstract}
Let $F_n:X_1\longrightarrow X_2$ be a sequence of (multivalued) meromorphic
maps between compact K\"ahler manifolds. We study 
the asymptotic distribution of 
preimages of points by $F_n$ and the asymptotic 
distribution of fixed points for multivalued self-maps
of a compact Riemann surface.
\par
Let $(Z_n)$ be a sequence of holomorphic images of $\P^s$ in a projective
manifold. We prove that the currents, 
defined by integration on $Z_n$,
properly normalized, converge to weakly laminar currents. We also 
show that the Green currents,
of suitable bidimensions, associated to a regular polynomial automorphism, are
(weakly) laminar. 
\end{abstract}
\section{Introduction}
L'une des motivations de notre travail est le probl\`eme dynamique 
suivant. 
Soit $f:X\longrightarrow X$ une application rationnelle ou birationnelle
sur une vari\'et\'e projective $X$ de dimension $k$. 
Soient $H_+$, $H_-$ des sous-vari\'et\'es
projectives de $X$, de dimensions respectives $s$, $k-s$ 
et de degr\'es fix\'es. 
Il s'agit de trouver des 
conditions suffisantes pour que la suite des mesures $\mu_n$, 
\'equidistribu\'ees aux 
points de $f^{-n}(H_+)\cap f^n(H_-)$, tende vers une mesure invariante $\mu$.
Nous pensons que la
mesure $\mu$ ne d\'ependra pas de $(H_+,H_-)$ g\'en\'erique
et sera un objet dynamique int\'eressant.
\par
Notons $Y$ l'espace (de param\`etres) des $(H_+,H_-)$. On peut construire
les applications m\'eromorphes multivalu\'ees 
$F_n:X\longrightarrow Y$  
telles que 
l'image r\'eciproque de $(H_+,H_-)$ par $F_n$ soit 
\'egale \`a l'ensemble $f^{-n}(H_+)\cap f^n(H_-)$.
On peut alors ramener le probl\`eme \`a l'\'etude 
de la distribution asymptotique des pr\'eimages des points par 
$F_n$ \cite{DinhSibony3}.
\par
On peut \'egalement approcher ce probl\`eme d'une autre mani\`ere
\cite{BedfordLyubichSmillie, Dujardin}. On 
cherche \`a d\'emontrer que les suites des courants d'int\'egration sur 
les vari\'et\'es
$f^{-n}(H_+)$ et $f^n(H_-)$, proprement normalis\'es, convergent vers 
des courants invariants $T_+$ et $T_-$. 
On attend que ces courants limites conservent des 
structures analytiques: ils sont constitu\'es par des familles de
vari\'et\'es complexes.
On attend aussi que la mesure invariante $T_+\wedge T_-$, 
co\"{\i}ncide avec la mesure obtenue comme intersection 
g\'eom\'etrique des vari\'et\'es qui constituent $T_+$ et $T_-$. 
On aura alors 
$\lim\mu_n = T_+\wedge T_-$.
\par
Ces deux diff\'erentes approches ont donn\'e 
\cite{BedfordLyubichSmillie,DinhSibony3, Dujardin2}
des r\'eponses partielles au probl\`eme ci-dessus dans 
le cas des automorphismes 
polynomiaux r\'eguliers au sens de Sibony. Notons qu'en dimension 2 
les automorphismes polynomiaux r\'eguliers sont ceux 
du type H\'enon (\voir paragraphe 5 pour 
la d\'efinition). 
\par
\
\par
Dans la premi\`ere partie de l'article, nous \'etudions la distribution 
des pr\'eimages d'une suite d'applications.
Soient $(X_1,\omega_1)$ et $(X_2,\omega_2)$ des vari\'et\'es 
k\"ahl\'eriennes compactes de dimensions respectives $k_1$ et $k_2$. 
Soit $F_n:X_1\longrightarrow X_2$ 
des applications m\'eromorphes (multivalu\'ees). 
Notons $d_n$ le degr\'e topologique et $\lambda_n$ 
le degr\'e interm\'ediaire d'ordre $k_2-1$ de $F_n$ (\voir paragraphe 2 pour 
les d\'efinitions).
On suppose que la s\'erie 
$\sum \lambda_nd_n^{-1}$ converge. Sodin, Russakovskii, Shiffman
\cite{RussakovskiiShiffman} ont montr\'e, pour le cas des 
applications rationnelles
entre espaces projectifs, que les pr\'eimages de $F_n$ sont \'equidistribu\'ees
lorsque $n$ tend vers l'infini. Plus pr\'ecis\'ement, 
$d_n^{-1}F_n^*(\delta_z)-d_n^{-1}F_n^*(\delta_{z'})$ tend  
faiblement vers 0 pour 
$z,z'\in X_2$ hors d'un ensemble pluripolaire $\E$, o\`u 
on a not\'e $\delta_z$ la masse
de Dirac en $z$. Le cas g\'en\'eral de ce r\'esultat a \'et\'e 
d\'emontr\'e dans \cite{DinhSibony3}. Si $(F_n)$ est la suite des it\'er\'es
d'une application  
m\'eromorphe (multivalu\'ee) 
$f$ dont le degr\'e topologique est plus grand 
que les autres degr\'es dynamiques, $\E$ est contenu dans une r\'eunion 
d\'enombrable d'ensembles analytiques, \voir Lyubich \cite{Lyubich}, 
Briend-Duval \cite{BriendDuval2, BriendDuval1}, Guedj \cite{Guedj} et 
\cite{DinhSibony2, Dinh,DinhSibony3}. De plus, $d_n^{-1}F_n^*(\delta_z)$ 
converge vers la mesure 
d'\'equilibre $\mu$ de $f$. Cette mesure refl\`ete aussi
la distribution des points fixes r\'epulsifs de $F_n$ ($=f^n$).
\par
Nous cherchons des conditions plus g\'en\'erales
sur la suite $(F_n)$ pour que $\E$ soit contenu 
dans une r\'eunion d\'enombrable 
d'ensembles analytiques. Nous montrons que c'est le cas si 
la suite des courants
postcritiques $S_n$ de $F_n$ 
($=1/d_n$ fois le courant d'int\'egration sur les 
valeurs critiques de $F_n$) converge vers un courant $S_\infty$. L'ensemble
$\E$ est alors
essentiellement contenu dans $\{\nu(S_\infty,z)>0\}$ o\`u $\nu(S_\infty,z)$
d\'esigne le nombre de Lelong de $S_\infty$ en $z$. Un th\'eor\`eme
de Siu \cite{Siu} implique que 
ce dernier ensemble est une r\'eunion finie ou d\'enombrable
d'ensembles analytiques.
\par
Nous \'etudions aussi la distribution des points fixes r\'epulsifs dans le cas
o\`u $X_1$ et $X_2$ sont \'egales \`a une surface de Riemann compacte $X$.
Supposons de plus que, pour tout $z\not\in\E$, 
$d_n^{-1}F_n^*(\delta_z)$ tend faiblement vers une 
mesure $\mu$ qui est singuli\`ere par rapport \`a $S_\infty$.
Nous montrons que la mesure, \'equidistribu\'ee aux 
points fixes r\'epulsifs de $F_n$, tend faiblement vers $\mu$. 
Ce r\'esultat \'eclaire le lien entre les ensembles postcritiques et
l'\'equidistribution des points fixes r\'epulsifs.
En donnant un exemple, nous 
montrons que l'hypoth\`ese sur le rapport entre $\mu$ et $S_\infty$ est 
n\'ecessaire. 
\par
La d\'emonstration des r\'esultats ci-dessus 
utilise une m\'ethode d\'evelopp\'ee dans 
Lyubich \cite{Lyubich},
Briend-Duval \cite{BriendDuval2} et \cite{DinhSibony2, Dinh}. Il s'agit de 
construire, pour des petites boules centr\'ees en un point 
g\'en\'erique $x_2\in X_2$,
un bon nombre de $F_n$-branches inverses dont on contr\^ole la taille. 
Plus pr\'ecis\'ement, 
nous construisons des applications inverses locales de $F_n$
qui sont d\'efinies sur les petites boules et qui admettent, comme images, 
des ensembles de petit diam\`etre.
L'\'equidistribution des pr\'eimages s'en d\'ecoule.
Pour construire les points fixes r\'epulsifs de $F_n:X\longrightarrow X$, nous 
utilisons une id\'ee de Lyubich \cite{Lyubich}. Si $U$ est un ouvert connexe
et simplement connexe tel que $\mu(U)>1-\epsilon/4$ et $S_\infty(U)<<\epsilon$, 
nous construisons $(1-\epsilon)d_n$ $F_n$-branches 
inverses de $U$ dont les images 
sont strictement contenues dans $U$. Une telle $F_n$-branche inverse 
est contractible pour
la m\'etrique de Kobayashi sur $U$ et, par cons\'equent, cr\'ee 
un point fixe r\'epulsif pour $F_n$.
\par
\
\par
Dans la deuxi\`eme partie de l'article, nous \'etudions la laminarit\'e 
de certains courants positifs ferm\'es. 
Soit $(Z_n)$ une suite d'images de $\P^s$ dans une vari\'et\'e projective 
$X$ de dimension $k$. Supposons 
que la suite des courants 
d'int\'egration sur $Z_n$, proprement 
normalis\'es, converge vers un courant $T$.
Nous montrons que 
$T$ est faiblement laminaire 
et est laminaire si $s=k-1$ et si les singularit\'es de $Z_n$ 
sont raisonnables. 
Rappelons la notion de laminarit\'e introduite par Bedford-Lyubich-Smillie
\cite{BedfordLyubichSmillie, Dujardin,deThelin}. On dit qu'un 
courant positif $T$ de bidimension 
$(s,s)$ est {\it faiblement laminaire}
s'il est localement \'ecrit comme une int\'egrale de courants d'int\'egration 
sur des vari\'et\'es complexes de dimension $s$. Si ces vari\'et\'es sont des
graphes disjoints, $T$ est dit {\it laminaire}
(\voir paragraphe 5 pour les d\'etails). 
\par
Nous utilisons aussi la construction de
branches inverses pour prouver ce r\'esultat.  
Plus pr\'ecis\'ement, nous consid\'erons 
des projections $F_n$ de $Z_n$ sur un espace projectif $\P^s$. Les 
$F_n$-branches 
inverses des ouverts de $\P^s$, avec taille contr\^ol\'ee, forment des familles
normales de vari\'et\'es complexes dans $X$. 
En passant \`a la limite,
on obtient des vari\'et\'es complexes qui constituent le courant $T$.
\par 
Comme application, nous montrons que les courants de Green, de certaines 
bidimensions,
d'un automorphisme polynomial 
r\'egulier, est laminaire ou faiblement laminaire.
Nos r\'esultats g\'en\'eralisent des th\'eor\`emes
de Bedford-Lyubich-Smillie \cite{BedfordLyubichSmillie} et de Dujardin 
\cite{Dujardin} qui ont \'etudi\'e le cas des courants 
de bidimension $(1,1)$. 
La laminarit\'e de courants de Green est aussi valable pour les automorphismes
d'une vari\'et\'e projective quelconque.
Notons que dans le cas des automorphismes du type H\'enon, la 
laminarit\'e de courants de Green a permis de d\'emontrer de nombreuses 
propri\'et\'es
dynamiques importantes
(\voir aussi \cite{Cantat, DillerFavre, McMullen, Dujardin3} 
pour le cas des applications bim\'eromorphes sur les surfaces). 
Nous voulons enfin signaler que dans un travail r\'ecent \cite{deThelin}, 
de Th\'elin a montr\'e que les limites de surfaces  de Riemann ouvertes,
de genres contr\^ol\'es, sont aussi des courants laminaires.
\section{Quelques d\'efinitions}
Nous allons introduire dans ce paragraphe quelques notions sur
les (semi)-transformations m\'eromorphes entre vari\'et\'es k\"ahl\'eriennes 
compactes.
La notion de transformations m\'eromorphes 
que nous utilisons ici correspond aux transformations de codimension 0 
de \cite{DinhSibony3}. 
\par
Le lecteur trouvera aussi dans ce paragraphe les principales notations
utilis\'ees dans tout l'article. Les 
vari\'et\'es k\"ahl\'eriennes compactes 
$(X,\omega)$, $(X_1,\omega_1)$ et 
$(X_2,\omega_2)$ sont de  
dimensions respectives $k$, $k_1$ et $k_2$. 
Leurs formes de K\"ahler sont normalis\'ees par 
$\int_X\omega^k=\int_{X_1}\omega_1^{k_1}=
\int_{X_2}{\omega_2}^{k_2}=1$.
$\Pi_1$ et $\Pi_2$ sont les projections canoniques de $X_1\times X_2$ ou de 
$X\times X$ sur le premier et le second facteur. 
Les espaces projectifs complexes sont muni de la forme de Fubini-Study 
normalis\'ee
$\omega_\FS$. Les notations $\delta_z$, 
$\aire(.)$, $\diam(.)$, $\trace(.)$ d\'esignent la masse de Dirac, l'aire,
le diam\`etre et la trace. La notation $\|\ \|$ d\'esigne la masse d'un courant,
la norme d'un vecteur ou d'un op\'erateur lin\'eaire.
\\
\
\\
{\bf 2.1. (Semi)-transformations m\'eromorphes} 
\\
\ 
\\
On appelle {\it $p$-cha\^{\i}ne holomorphe} de $X_1\times X_2$ 
toute combinaison
lin\'eaire $\Gamma:=\Gamma^1+\cdots +  \Gamma^m$ o\`u les $\Gamma^i$
sont des sous-ensembles analytiques irr\'eductibles de dimension $p$ de 
$X_1\times X_2$. On ne suppose pas 
que les $\Gamma^i$ sont distincts. Notons $|\Gamma|:=\cup \Gamma_i$ le support
de $\Gamma$ et $[\Gamma]:=\sum[\Gamma^i]$ 
le courant d'int\'egration sur $\Gamma$. {\it La multiplicit\'e} 
$\mult(\Gamma,z)$ de $\Gamma$
en un point $z$ est le nombre de $\Gamma^i$ qui contiennent $z$. 
\par
On appelle {\it semi-transformation m\'eromorphe} (STM) de $X_1$ dans $X_2$ 
la donn\'ee d'une $k_2$-cha\^{\i}ne holomorphe
$\Gamma$ de $X_1\times X_2$ telle que la restriction de $\Pi_2$
\`a chaque composante irr\'eductible de $\Gamma$ soit surjective.
On identifie cette STM \`a ``l'application multivalu\'ee''
$F:=\Pi_2\circ(\Pi_{1|\Gamma})^{-1}$ et on dira que $\Gamma$ est 
{\it le graphe} de
$F$. 
Lorsque la restriction de $\Pi_1$
\`a chaque composante irr\'eductible de $\Gamma$ est surjective, $F$ est 
une {\it transformation m\'eromorphe} (TM). 
La TM $F^{-1}:X_2\longrightarrow X_1$ de graphe $\Gamma$ est appel\'ee 
{\it TM adjointe} de $F$.
Une TM entre vari\'et\'es 
de m\^eme dimension est appel\'ee {\it correspondance m\'eromorphe} (CM).
Si $F$ est une CM telle que le degr\'e de $\Pi_{1|\Gamma}$ soit
\'egal \`a 1, $F$ est une application m\'eromorphe surjective de $X_1$ dans $X_2$.
\par
Soit $F:X_1\longrightarrow X_2$ une STM. 
Notons $\I$ l'ensemble des points $x_2\in X_2$ tels que 
$\dim \Pi_2^{-1}(x_2)\cap \Gamma\geq 1$. C'est un sous-ensemble analytique 
de codimension au moins 2 de $X_2$. Il est appel\'e {\it deuxi\`eme ensemble
d'ind\'etermination} de $F$. 
On appelle {\it degr\'e topologique} de $F$ le nombre $d_t$ de 
points de $\Pi_2^{-1}(x_2)\cap\Gamma$, compt\'e avec multiplicit\'e, pour
$x_2\in X_2\setminus \I$. Ce nombre ne d\'epend pas de $x_2$. Il est 
aussi \'egal au nombre de points de $F^{-1}(x_2)$ compt\'e avec 
multiplicit\'e o\`u $F^{-1}:=\Pi_1\circ(\Pi_{2|\Gamma})^{-1}$.
\par
Posons $F^*:=(\Pi_1)_*(\Pi_{2|\Gamma})^*$ et 
$F_*:=(\Pi_2)_*(\Pi_{1|\Gamma})^*$. L'op\'erateur $F^*$ (resp. $F_*$) 
agit sur l'espace des formes lisses sur $X_2$ 
(resp. sur $X_1$) \`a valeurs dans l'espace des courants 
sur $X_1$ (resp. sur $X_2$). L'op\'erateur $F^*$ agit 
aussi sur les mesures positives $\nu$ 
qui ne chargent pas 
$\I$. La masse de $F^*(\nu)$ est $d_t$ fois celle de $\nu$. On appelle 
{\it degr\'e interm\'ediaire d'ordre $p$} de $F$ le nombre
$$\delta_p(F):=\int [\Gamma]\wedge \Pi_1^*(\omega_1^{k_2-p})\wedge 
\Pi_2^*(\omega_2^p),\ \ 0\leq p\leq k_2.$$
On a 
$$\delta_p(F):=\int_{X_2}F_*(\omega_1^{k_2-p})\wedge \omega_2^p.$$
Observons que les degr\'es interm\'ediaires se calculent cohomologiquement.
\par
Pour chaque point $z\in\Gamma^i$, 
notons $\Gamma^i_j(z)$ les germes d'ensemble analytique irr\'eductible de 
$\Gamma^i$ en $z$. Soit $m_{i,j}(z)$ le degr\'e topologique de
$\Pi_{2|\Gamma^i_j(z)}$. Posons $n_i(z):=\sum (m_{i,j}(z)-1)$.
L'ensemble $\{n_i\geq 1\}$ est
une hypersurface de $\Gamma^i$. Soit 
${\rm C}^i$ la  $(k_2-1)$-cha\^{\i}ne holomorphe port\'ee par cette hypersurface
dont les multiplicit\'es sont donn\'ees par la fonction $n_i$.
Posons ${\rm C}:=\sum {\rm C}^i$.
On appelle {\it courant postcritique} de $F$ le courant 
$S:=(\Pi_2)_*[{\rm C}]$.  
\par
Soit $K_0$ un sous-ensemble connexe de $X_2$. On appelle {\it branche 
inverse} de $K_0$ toute suite $\B$ 
$$K_{-1},(\widehat K_{-1},i),K_0$$
avec $\widehat K_{-1}$ un sous-ensemble connexe de $\Gamma^i$, 
$K_{-1}:=\Pi_1(\widehat K_{-1})$ 
telle que $\Pi_2$ d\'efinisse une bijection entre 
$\widehat K_{-1}$ et $K_0$. 
Si $K_0$ n'est pas un ouvert de $X_2$, on exige que 
l'inverse $\tau^\B$ de l'application $\Pi_2: \widehat K_{-1}\longrightarrow 
K_0$ se prolonge en application holomorphe 
d'un voisinage de $K_0$ dans $\Gamma^i$. 
Si $K_0$ est un ouvert de $X_2$, on exige bien s\^ur que $\tau^\B$ est holomorphe. 
On dira que $\tau^\B$ (resp. $F^\B:=\Pi_1\circ \tau^\B$) est {\it 
l'application 
semi-inverse} (resp. {\it inverse}) associ\'ee \`a la branche $\B$.
On dira que $\B$ est {\it de taille} $\delta$ si $\diam (K_{-1})\leq
\delta$. L'ensemble $K_{-1}$ est appel\'e {\it l'image} de $\B$. 
Consid\'erons une branche inverse $\B_0$ d'un point $x_0\in X_2$:
$$x_{-1}, (\widehat x_{-1},j), x_0.$$
On dira que $\B$ est {\it accroch\'ee} \`a $\B_0$ si 
$x_0\in K_0$, $j=i$ et $\widehat x_{-1}\in\widehat K_{-1}$. Etant donn\'ee
une branche inverse $\B_0$ 
de $x_0$, par prolongement analytique, il y a au plus
une branche inverse de $K_0$ qui est accroch\'ee \`a $\B_0$. En particulier, 
$K_0$ admet au plus $d_t$ branches inverses. Observons que ${\rm C}$ et $S$ 
sont des obstructions pour construire les branches inverses des ensembles 
simplement connexes. Par contre, les points singuliers de $\Gamma^i$, dont 
les composantes irr\'eductibles locales sont lisses, ne le sont pas.
\par
\ 
\par
Consid\'erons maintenant le cas o\`u $X_1=X_2=X$ et o\`u $f$ est
une CM de degr\'e topologique $d_t$ de $X$ dans $X$. 
On d\'efinit {\it l'it\'er\'e 
d'ordre $n$} de $f$ par 
$f^n:=f\circ \cdots \circ f$ ($n$ fois) o\`u la composition consid\'er\'ee est 
celle des 
applications multivalu\'ees. On appelle {\it degr\'e dynamique d'ordre
$p$} de $f$ le nombre
$$d_p:=\limsup_{n\rightarrow\infty} [\delta_p(f^n)]^{1/n}.$$
Les degr\'es dynamiques de $f$ ne d\'ependent pas de la forme de K\"ahler fix\'ee
pour $X$ et on a $d_k=d_t$. 
\par
Les points fixes de $f$ sont d\'etermin\'es par l'intersection de 
son graphe $\Gamma$
avec la diagonale $\Delta$ de $X\times X$. 
Soit $\B$ une branche inverse d'un point $x\in X$ et soit $f^\B$ l'application
inverse associ\'ee. Lorsque $f^\B(x)=x$ on dit que $\B$ {\it d\'efinit  
un point fixe r\'egulier} de $f$. Si, de plus, les valeurs propres 
de $f^\B$ en $x$ sont de module strictement inf\'erieur \`a 1, on dit que ce
point fixe est {\it r\'epulsif}. Quand il n'y a pas de confusion, on note
simplement $x$ pour le point fixe. 
Les {\it points p\'eriodiques r\'eguliers (r\'epulsifs) 
d'ordre $n$} de $f$ sont les points fixes r\'eguliers (r\'epulsifs) de $f^n$. 
\\
\
\\
{\bf 2.2. Suites de (semi)-transformations m\'eromorphes}
\\
\ 
\\ 
Consid\'erons une suite de STM 
$F_n:X_1\longrightarrow X_2$. Soient
$d_n$ le degr\'e
topologique et
$\lambda_n$ le degr\'e interm\'ediaire d'ordre $k_2-1$ de $F_n$. Notons $S_n$
le courant postcritique et $\I_n$ le deuxi\`eme 
ensemble d'ind\'etermination de $F_n$. 
\par
Nous allons consid\'erer les  suites $(F_n)$ v\'erifiant une ou 
plusieurs des propri\'eti\'es suivantes. 
\begin{enumerate}
\item[(H1) ] La suite $\lambda_nd_n^{-1}$ tend vers $0$;
\item[(H1')] La s\'erie $\sum_{n\geq 0}\lambda_nd_n^{-1}$ est convergente;
\item[(H2) ] La suite de mesures 
$d_n^{-1} F_n^*({\omega_2}^{k_2})$ tend faiblement vers une
mesure $\mu$; on dira qu'alors $\mu$ est la {\it mesure d'\'equilibre}
de la suite $(F_n)$;
\item[(H3) ] La suite $(S_n)$ tend faiblement vers un courant 
$S_\infty$;
on dira qu'alors $S_\infty$ est 
le {\it courant postcritique} de la suite $(F_n)$;
\end{enumerate}
Notons $\E_1$ 
l'ensemble des points $x_2\in X_2$ qui 
appartiennent \`a $\I_n$ pour une infinit\'e de $n$ et $\E_2$ l'ensemble 
des points $x_2\in X_2\setminus \E_1$ tels que 
la mesure
$\mu_n^{x_2}:=d_n^{-1}F_n^*(\delta_{x_2})$ 
ne tende pas vers $\mu$ quand $n\rightarrow\infty$.
L'ensemble $\E_1$ (resp. $\E_2$) est 
appel\'e {\it premier} (resp. {\it deuxi\`eme}) 
{\it ensemble exceptionnel} de $(F_n)$. On dira aussi que $\E:=\E_1\cup \E_2$ 
est {\it l'ensemble exceptionnel} de $(F_n)$. 
Posons $\E_1^*:=\cup_{n\geq 0} \I_n$. C'est une r\'eunion finie ou 
d\'enombrable
d'ensembles analytiques de codimension au moins 2 de $X_2$. 
On a $\E_1\subset \E_1^*$.  
Dans \cite{DinhSibony3}, en g\'en\'eralisant des r\'esultats de Sodin, 
Russakovskii-Shiffman \cite{RussakovskiiShiffman}, nous
avons montr\'e, pour toute suite de TM 
v\'erifiant (H1') et (H2), que l'ensemble exceptionnel est pluripolaire.
\section{Ensemble exceptionnel}
Dans ce paragraphe, nous d\'emontrons que si une suite 
de TM v\'erifie (H1), (H2), (H3), alors son ensemble 
exceptionnel 
est contenu dans une r\'eunion 
finie ou d\'enombrable d'ensembles analytiques. Nous avons le
r\'esultat suivant. 
\begin{theoreme} Soit $(F_n)$ une suite  
de TM entre vari\'et\'es k\"ahl\'eriennes compactes
$(X_1,\omega_1)$ et $(X_2,\omega_2)$. 
Supposons que $(F_n)$ v\'erifie (H1), (H2), (H3). 
Soient $S_\infty$ le courant postcritique et 
$\E_2$ le deuxi\`eme ensemble exceptionnel 
de $(F_n)$. Alors $\E_2$ est contenu dans l'ensemble
$\E^*_2:=\{\nu(S_\infty,x_2)>0\}$ 
o\`u $\nu(S_\infty,x_2)$ d\'esigne 
le nombre de Lelong de $S_\infty$ en $x_2$. 
\end{theoreme}
\par
Observons que pour toute suite $(F_n)$ v\'erifiant (H1), (H3), 
on peut extraire des sous-suites v\'erifiant (H2). Par cons\'equent,
pour une telle suite, $\mu_n^{x_2}-\mu_n^{x_2'}\rightarrow 0$ lorsque
$x_2,x_2'\not\in\E_1\cup \E_2^*$.
\par
Soit $f$ une CM d'une vari\'et\'e $X$ comme au paragraphe 2.1. 
Si son degr\'e 
topologique est strictement plus grand que 
son degr\'e dynamique d'ordre $k_2-1$,
on peut appliquer le th\'eor\`eme 3.1 \`a la suite des it\'er\'es de $f$.
Dans ce cas, la suite $(S_n)$ est croissante, born\'ee en masse. Donc elle
converge vers un courant $S_\infty$.
On obtient alors des r\'esultats de
Lyubich \cite{Lyubich}, Freire-Lopes-Ma\~n\'e \cite{FreireLopesMane}, 
Briend-Duval \cite{BriendDuval2} et \cite{DinhSibony2,Guedj,Dinh, DinhSibony3}.
\par
Dans la d\'emonstration du th\'eor\`eme 3.1, 
nous construisons,  pour chaque $F_n$, 
un bon nombre de branches
inverses sur des petites boules centr\'ees en un point 
$x_2\in X_2\setminus \E_1\cup \E_2^*$. 
La taille de ces $F_n$-branches inverses 
tend vers
0 quand $n$ tend vers l'infini et le th\'eor\`eme 3.1 s'en d\'ecoule. 
Nous commen\c cons par la construction de branches inverses
pour une STM g\'en\'erale 
$F:X_1\longrightarrow X_2$. 
Notons $\Gamma=\sum\Gamma^i$, 
$d$ et $\lambda$ le graphe, 
le degr\'e topologique et le degr\'e interm\'ediaire
d'ordre $k_2-1$ de $F$. Notons $\I$ le deuxi\`eme ensemble d'ind\'etermination 
et $S$ le courant postcritique de $F$. Posons 
$\Omega:=d^{-1}F_*(\omega_1)$. Fixons un point $x_2\in X_2\setminus \I$. 
Pour simplifier les notations et les calculs, nous allons 
consid\'erer une carte de $X_2$ contenant $x_2$.
Dans cette carte, 
les boules et la masse de courant seront 
d\'efinies ou mesur\'ees avec la m\'etrique euclidienne. Ceci ne changera pas 
le r\'esultat. 
\par
Soit $\delta_0>0$ une constante assez petite, 
que nous pr\'eciserons plus loin.
Elle ne d\'epend que de la g\'eom\'etrie de
$(X_1,\omega_1)$. 
Soit $\Ball_{r_0}$ la boule 
de centre $x_2$ et de rayon $r_0$. Fixons des constantes  
$\nu\geq 0$ et $\delta\geq 0$ telles que
la masse de $S$ et de $\Omega$ 
dans  $\Ball_{r_0}$
soient major\'ees par $r_0^{2k_2-2}\nu$ et par $r_0^{2k_2-2}\delta$.
La constante $A$, que nous utilisons dans la proposition suivante,  
sera donn\'ee au lemme 3.3.
\begin{proposition} 
Soit $F:X_1\longrightarrow X_2$ une STM comme ci-dessus.
Soit $\epsilon$ une constante, $0<\epsilon<1$. 
Supposons que $\nu<A^{-1}\epsilon^2/8$ et $\delta<\delta_0$.
Alors il existe $r$, $0<r\leq r_0$, 
ind\'ependant de $F$ 
tel que la boule 
$\Ball_r$ admette au moins 
$(1-\epsilon)d$ branches inverses de taille $\delta^{1/2}$.
\end{proposition} 
{\bf D\'emonstration du th\'eor\`eme 3.1.} Soit 
$x_2\in X_2\setminus (\E_1\cup\E_2^*)$. On a $\nu(S_\infty,x_2)=0$. 
Il faut montrer que $\mu^{x_2}_n$ 
tend faiblement vers $\mu$. 
Quitte \`a extraire une sous-suite, on peut supposer que $(F_n)$
v\'erifie (H1'), 
$x_2\not\in \E_1^*$ et que $\mu^{x_2}_n$ converge vers une mesure 
$\mu^{x_2}$. 
\par
Fixons des constantes 
$\epsilon$ et $\nu$ v\'erifiant l'hypoth\`ese de la proposition 3.2.
D'apr\`es 
la d\'efinition du nombre de Lelong
\cite{Skoda} et puisque $S_n\rightarrow S_\infty$, on
peut choisir 
$r_0>0$ assez petit tel que, pour $n$ assez grand,
la masse de $S_n$ dans $\Ball_{r_0}$ soit major\'ee par 
$r_0^{2k_2-2}\nu$.
Posons $\Omega_n:=d_n^{-1}(F_n)_*\omega_1$.
La masse de ce courant est d'ordre 
$\lambda_nd_n^{-1}$.
Puisque $(F_n)$ v\'erifie (H1), on peut appliquer la proposition 3.2 
pour $n$ assez grand.
Notons $\B_{n,j}$ les 
$F_n$-branches inverses de $\Ball_r$ 
fournies par la proposition 3.2 et $F^{\B_{n,j}}$ 
les applications inverses associ\'ees. 
La condition (H1)
implique que la taille de $\B_{n,j}$ 
tend vers 0 quand $n$ tend vers l'infini.
\par
D'apr\`es \cite{DinhSibony3}, puisque $(F_n)$ v\'erifie (H1') et (H2), 
$\E$ est pluripolaire.
On peut donc choisir un point $x_2'\in \Ball_r$ tel que 
$\mu_n^{x_2'}\rightarrow \mu$.
Posons 
$$\tilde\mu^{x_2}_n:=d_n^{-1}\sum_j 
(F^{\B_{n,j}})_*(\delta_{x_2})\ \mbox{ et }\  
\tilde\mu^{x_2'}_n:=d_n^{-1}\sum_j (F^{\B_{n,j}})_*(\delta_{x_2'}).$$ 
Quitte \`a extraire des
sous-suites, on peut supposer que $(\tilde\mu^{x_2}_n)$
et $(\tilde\mu^{x_2'}_n)$ sont convergentes. 
Le fait que la taille de $\B_{n,j}$ tend vers 0
implique que ces deux suites ont une m\^eme limite.
D'autre part,
pour chaque $n$ assez grand, le nombre des branches $\B_{n,j}$ est au moins 
\'egal \`a $(1-\epsilon)d_n$. Les masses des 
mesures positives $(\mu^{x_2}-\lim  \tilde\mu^{x_2}_n)$ et 
$(\mu-\lim  \tilde\mu^{x_2'}_n)$
sont donc au plus \'egales \`a $\epsilon$. 
On en d\'eduit que la masse de $\mu^{x_2}-\mu$ est major\'ee par 
$2\epsilon$ pour tout $\epsilon>0$.
D'o\`u $\mu^{x_2}=\mu$. 
\par
\hfill $\square$
\\
\par
Dans la suite, nous d\'emontrons la proposition 3.2.
Pour simplifier la construction des branches inverses, on suppose que $|\Gamma|$
est localement irr\'eductible en tout point. Pour le cas contraire, il suffit 
d'utiliser une application $\tau:\tilde\Gamma\longrightarrow\Gamma$ afin de
s\'eparer les composantes irr\'eductibles locales de $\Gamma$ et on remplace
$\Pi_1$, $\Pi_2$ par $\Pi_1\circ\tau$ et $\Pi_2\circ\tau$.
\par
Soit $\G$ la famille des droites
complexes passant par $x_2$. On munit $\G$ de la structure 
de vari\'et\'e complexe naturelle. 
Puisque $\G\simeq \P^{k_2-1}$, on munit $\G$ de 
la mesure de probabilit\'e invariante
naturelle $\H^{2k_2-2}$. 
Notons $\Delta_\xi(r)$ le disque de rayon $r$ centr\'e en $x_2$ 
et contenu dans la droite $\Delta_\xi$, $\xi\in\G$. Notons aussi
$[\Delta_\xi(r)]$
le courant d'int\'egration sur $\Delta_\xi(r)$. 
\begin{lemme} Soit $T$ un courant 
positif ferm\'e de bidegr\'e $(1,1)$, de masse 
$r^{2k_2-2}M$ dans la boule $\Ball_r$. Notons 
$\G(T,c)$ la famille des droites $\Delta_\xi$ telles que la masse de 
la mesure $T\cap[\Delta_\xi(r)]$ soit major\'ee par $c$. Alors
il existe une constante $A>0$
ind\'ependante de $\epsilon$, 
$T$, $M$, $r$ telle que $\H^{2k_2-2}(\G(T,A\epsilon^{-1} M))\geq 1-
\epsilon/4$.
\end{lemme}
\begin{preuve} Par homoth\'etie, on peut supposer que $r=M=1$. Soient 
$\pi: \widehat \Ball_2\longrightarrow \Ball_2$ l'\'eclatement de $\Ball_2$ 
en $x_2$ et $\widehat \Ball_1:=\pi^{-1}\Ball_1$. 
Fixons une m\'etrique k\" ahl\'erienne sur $\widehat \Ball_2$.
Dans $\Ball_1$, on peut \'ecrire $T=\ddc\varphi$ o\`u $\varphi$ est
une fonction p.s.h. D\'efinissons $\pi^*(T):=\ddc\varphi\circ \pi$. Ce
courant ne d\'epend pas du choix du potentiel $\varphi$ 
et l'op\'erateur $\pi^*$ 
est continu sur l'ensemble des courants $T$ \cite{Meo}. Par cons\'equent,
la masse de $\pi^*(T)$ est born\'ee par $c$ o\`u $c>0$
est une constante ind\'ependante 
de $T$.
\par
Notons $\widehat\Delta_\xi(1)$ 
l'adh\'erence de $\pi^{-1}(\Delta_\xi(1)\setminus\{x_2\})$. La masse de 
$T\cap[\Delta_\xi(1)]$ est \'egale \`a celle de 
$\pi^*(T)\cap[\widehat\Delta_\xi(1)]$. 
Puisque la fibration de $\widehat\Ball_1$, form\'ee par les disques 
$\widehat\Delta_\xi(1)$, n'a pas de singularit\'e, 
on peut appliquer la th\'eorie de 
tranchage classique pour conclure \cite[4.3.2]{Federer}.
\end{preuve}
\par
D'apr\`es le lemme 3.3, 
il existe une famille $\G'\subset \G$
avec $\H^{2k_2-2}(\G')\geq 1-\epsilon/2$ telle que pour tout
$\xi\in\G'$ la masse de la mesure 
$\Omega\cap[\Delta_\xi(r_0)]$ 
soit major\'ee par $A\epsilon^{-1}\delta$ et celle de  
$S\cap[\Delta_\xi(r_0)]$ soit major\'ee par 
$A\epsilon^{-1}\nu<\epsilon/8$. 
Pour construire les branches inverses 
de $\Delta_\xi(r)$, on a besoin du lemme suivant.
\begin{lemme} Soient $\Sigma$ une 
surface de Riemann \`a bord lisse et $\Delta_1$ le disque unit\'e. 
Soit $\pi:\Sigma\longrightarrow \Delta_1$ une application holomorphe 
lisse jusqu'au bord. Supposons que $\pi$ d\'efinisse un rev\^etement
de degr\'e $d$ sans points de ramification au bord. 
Soit $m$ le nombre de points de ramification de 
$\pi$ compt\'es avec multiplicit\'e.
Alors $\Sigma$ admet 
au moins $d-2m$ composantes 
connexes, irr\'eductibles, sur lesquelles $\pi$ est de degr\'e $1$.
\end{lemme}
\begin{preuve} Notons que tout point singulier de $\Sigma$ est un point
de ramification 
dans notre sens. Soient $\Sigma_1,\ldots,\Sigma_s$ 
les composantes connexes de 
$\Sigma$ et $\pi_i$ la restriction de $\pi$ \`a $\Sigma_i$.
Soient $\chi_i$ la caract\'eristique d'Euler de $\overline \Sigma_i$, 
$d_i$ le degr\'e de $\pi_i$ et $m_i$ le nombre de points de ramification de 
$\pi_i$ compt\'e avec multiplicit\'e. 
On a $\sum d_i=d$ et $\sum m_i= m$. 
D'apr\`es la formule de Riemann-Hurwitz \cite[p.105]{Kirwan}, on a
$\chi_i=d_i-m_i$ et donc $\sum\chi_i= d-m$. 
Puisque $\chi_i\leq 1$, on a $ \chi_i=1$  pour au moins 
$d-m$ indices $i$. On en d\'eduit que $\chi_i=1$ et 
$m_i=0$ pour au moins $d-2m$
indices $i$. Pour un tel indice $i$, on a $d_i=\chi_i+m_i=1$. 
\end{preuve}
\begin{lemme} Il existe $r_1>0$ ind\'ependant de $F$ tel que 
pour tout $\xi\in\G'$, $\Delta_\xi(r_1)$ 
admette au moins $(1-\epsilon/2)d$ branches inverses de taille 
$\frac{1}{2}\delta^{1/2}$.
\end{lemme}
\begin{preuve} Soit $\Sigma^i$ l'image r\'eciproque de $\Delta_\xi(r_0)$ par 
$\Pi_{2|\Gamma^i}$. Quitte \`a perturber l\'eg\`erement $r_0$, on peut supposer
que le bord de $\Sigma^i$ est lisse et ne contient pas de point
de ramification de $\Pi_{2|\overline\Sigma^i}$. 
Soit $m_i$ le nombre de points de ramification 
de $\Pi_{2|\Sigma^i}$. La somme $m:=\sum m_i$  
est \'egale \`a $d$ fois la masse de 
$S\cap[\Delta_\xi(r_0)]$. Donc $m\leq\epsilon d/8$. Le lemme 3.4, appliqu\'e
aux $\Sigma^i$, implique que  
$\Delta_\xi(r_0)$ admet au moins 
$(1-\epsilon/4)d$ branches inverses.
\par
L'aire totale de ces branches inverses est major\'ee par $d$ fois la masse de
la mesure $\Omega\cap[\Delta_\xi(r_0)]$. Elle est donc 
major\'ee par $A\epsilon^{-1}\delta d$.
On en d\'eduit que pour 
au moins $(1-\epsilon/2)d$ branches inverses de
$\Delta_\xi(r_0)$, l'aire est major\'ee par $4A\epsilon^{-2}\delta$.
La preuve est compl\'et\'ee par le lemme 3.6.
\end{preuve}
\begin{lemme} Soient $\delta_1>0$ une constante assez petite
et $\Delta_R\subset \C$ le disque de rayon $R$ centr\'e en $0$. Alors, si
$\tau:\Delta_1\longrightarrow X_1$ est une application holomorphe v\'erifiant 
$\aire(\tau(\Delta_1))\leq \delta_1$, on a 
\begin{enumerate}
\item Pour tout $c>0$, il existe $r$, $0<r<1$, ind\'ependant de
$\tau$, tel que  
$\diam(\tau(\Delta_r))\leq c\sqrt{\aire(\tau(\Delta_1))}$;
\item Pour tout $r$, $0<r<1$, il existe $c>0$, ind\'ependant de
$\tau$, tel que  
$\diam(\tau(\Delta_r))\leq c\sqrt{\aire(\tau(\Delta_1))}$;
\end{enumerate}
o\`u l'aire de $\tau(\Delta_1)$ est calcul\'ee en tenant compte 
les multiplicit\'es des points.
\end{lemme}
Ce lemme est d\^u \`a Briend-Duval. La preuve est donn\'ee 
pour le cas de l'espace projectif, elle est aussi 
valable pour le cas g\'en\'eral
\cite{BriendDuval2}. On peut 
\'egalement montrer ce lemme en consid\'erant 
une suite g\'en\'erale d'applications $\tau_n:\Delta_1
\longrightarrow X_1$ avec 
$\lim \aire(\tau_n(\Delta_1)) =0$. On montre que les 
valeurs adh\'erentes de cette suite
sont des applications constantes. Ceci permet 
d'appliquer la formule de Cauchy 
dans des cartes convenables pour conclure.
\par
\ 
\par
Nous recouvrons $X_1$ par un nombre fini de cartes biholomorphes
\`a la boule unit\'e de $\C^{k_1}$. 
Fixons la constante
$\delta_0$ telle que tout ensemble 
de diam\`etre $2\delta_0^{1/2}$ soit contenu
dans l'une des cartes. 
Dans la suite, nous allons travailler avec 
des sous-ensembles de $X_1$ de diam\`etre plus petit que $\delta_0^{1/2}$. 
Nous pouvons donc utiliser les m\'etriques euclidiennes sur 
ces cartes pour simplifier les notations et les calculs. 
\par  
Notons $\F$ la famille des branches inverses de $x_2$. Pour chaque
$\xi\in \G'$, notons $\F_{\xi}$ la famille des branches $\B\in \F$ qui 
accrochent une branche inverse de taille 
$\frac{1}{2}\delta^{1/2}$ 
de $\Delta_\xi(r_1)$. D'apr\`es le lemme 3.5, on a $\#\F_{\xi}\geq 
(1-\epsilon/2)d$. Notons $\G_{\B}$ la famille des $\xi\in\G'$ tels
que $\Delta_\xi(r_0)$ admette une branche inverse accroch\'ee \`a $\B$.
On a 
$$\sum_\B \H^{2k_2-2}(\G_{\B})\geq 
(1-\epsilon/2)d
\H^{2k_2-2}(\G')\geq (1-\epsilon/2)^2d.$$ 
Soit $m$ le nombre d'\'el\'ements 
$\B\in\F$ tels que $\H^{2k_2-2}(\G_{\B})\geq \epsilon/4$.
Puisque  $\H^{2k_2-2}(\G_{\B})\leq 1$ pour tout $\B$, on a
$$(1-\epsilon/2)^2d \leq \sum_\B \H^{2k_2-2}(\G_{\B}) \leq   
m+(d-m)\epsilon/4.$$
D'o\`u $m\geq (1-\epsilon)d$.
\par
Fixons une branche $\B$ parmi les $m$ branches qui v\' erifient 
$\H^{2k_2-2}(\G_{\B})\geq \epsilon/4$.
Soient $\tau^\B$ et $F^\B$ 
les applications semi-inverse et inverse  associ\'ees \`a $\B$. 
Ce sont des applications
holomorphes au voisinage de $x_2$.
Pour terminer la preuve de la proposition 3.2, 
il suffit d'appliquer le lemme 3.7 afin de prolonger 
$\tau^\B$ et $F^\B$ en applications holomorphes sur 
$\Ball(x_2,r)$ avec $r:=u r_1$. Les applications prolong\'ees d\'efinissent 
alors une branche inverse pour $\Ball(x_2,r)$ qui est de taille 
$\delta^{1/2}$. 
\par
\hfill $\square$
\begin{lemme}{\bf \cite{Alexander, SibonyWong}}
Soit $\G'$ une famille de droites passant par $x_2$ telle que  
$\H^{2k_2-2}(\G')\geq \epsilon/4$. Notons
$\Sigma$ l'intersection de ces droites avec
la boule $\Ball(x_2,r_1)$. 
Soit $\tau$ une application holomorphe d'un voisinage de $x_2$ dans
$\C^k$. Supposons que $\tau$ se prolonge holomorphiquement 
sur $\Delta_\xi\cap\Ball(x_2,r_1)$ pour tout $\xi\in\G'$.
Alors il existe $u>0$, d\'ependant de $\epsilon$, mais 
ind\'ependant de $\G'$ et de $\tau$,
tel que $\tau$ se prolonge 
holomorphiquement sur
$\Ball(x_2,ur_1)$. De plus, si on note encore $\tau$ les prolongements 
holomorphes, on a
$$\sup_{x_2'\in \Ball(x_2,ur_1)}\|\tau(x_2')-\tau(x_2)\| \leq
\sup_{x_2'\in \Sigma} \|\tau(x_2')-\tau(x_2)\|.$$
En particulier, on a $\diam \tau(\Ball(x_2,ur_1))\leq 2\diam \tau(\Sigma)$.
\end{lemme}
\par
Soient $H^i$ une $p$-cha\^{\i}ne holomorphe, $0\leq p\leq k_2-1$,
\`a support dans $\Gamma^i$.
Notons $H^i_*$ l'ensemble des points $z\in\Gamma^i$ 
tels que la multiplicit\'e de $H^i$ en $z$ est plus grande ou \'egale 
au nombre de composantes irr\'eductibles lisses de $\Gamma^i$ en $z$.
Posons $H:=\sum H^i$.
\begin{proposition} Sous l'hypoth\`ese de la proposition 3.2, il existe 
une constante $\nu_1>0$ ind\'ependante de $F$ telle que 
si la masse du courant
$d^{-1}(\Pi_2)_*[H]$ dans $\Ball_{r_0}$ est major\'ee
par $r_0^{2p} \nu_1\epsilon$, $\Ball_{r/2}$ admette au plus
$2\epsilon d$ branches inverses 
$$K_{-1}, (\widehat K_{-1},i), K_0$$
(avec $K_0=\Ball_{r/2}$)
qui v\'erifient $\widehat K_{-1}\cap H^i_*\not=\emptyset$.
\end{proposition}
\begin{preuve} 
D'apr\`es une in\'egalit\'e du type 
Jensen \cite{Skoda}, la masse de 
$d^{-1}(\Pi_2)_*[H]$ dans $\Ball_r$ est major\'ee
par $r^{2p} \nu_1\epsilon$. 
D'apr\`es une in\'egalit\'e de Lelong, lorsque $\nu_1$ est 
assez petit, tout sous-ensemble analytique de dimension $p$ 
de $\Ball_r$ qui rencontre $\Ball_{r/2}$, a un volume 
sup\'erieur \`a $\nu_1r^{2p}$.
Soit $\B$ une branche inverse de $\Ball_r$ donn\'ee par la Proposition 3.2.
Notons $\B'$ la branche inverse   
$$K_{-1}, (\widehat K_{-1},i), K_0$$
de $\Ball_{r/2}$ ($=K_0$) induite par $\B$. 
Si la branche inverse $\B'$ v\'erifie 
$H^i_\B:= \widehat K_{-1}\cap H^i_*\not=\emptyset$,
la masse de $(\Pi_2)_*[H^i_\B]$ est 
minor\'ee par $\nu_1r^{2p}$. Puisque la masse 
totale des $(\Pi_2)_*[H^i_\B]$ est minor\'ee par $dr^{2p}\nu_1\epsilon$,
il y a au plus $\epsilon d$ branches inverses 
$\B$ v\'erifiant cette propri\'et\'e.
La proposition en d\'ecoule.
\end{preuve}
\section{Distribution de points p\'eriodiques}
Consid\'erons une suite de CM $(F_n)_{n\geq 0}$ 
de $X$ dans lui-m\^eme. Supposons que $(F_n)$ v\'erifie (H2).
Notons $\FR_n$ la famille des points fixes r\'eguliers r\'epulsifs de $F_n$. 
Dans cette famille, on r\'ep\`ete chaque point un nombre de fois \'egal \`a
sa multiplicit\'e. Posons 
$$\widehat\mu_n:=\sum_{x\in \FR_n}\delta_x.$$
Nous cherchons des conditions suffisantes pour que 
$\widehat\mu_n$
tende faiblement vers la mesure d'\'equilibre $\mu$ de $(F_n)$. 
Nous donnons dans le cas de dimension 1 un crit\`ere 
satisfaisant. 
\begin{theoreme} Soit $X$ une surface de Riemann compacte lisse. Soit
$(F_n)$ une suite de correspondances de $X$ dans $X$. 
Supposons que $(F_n)$ v\'erifie (H1), (H2), (H3) et 
$$\mbox{\it (H4)\ \ \ \ la mesure $S_\infty$ est 
singuli\`ere par rapport \`a } \mu.$$
Alors $\widehat\mu_n$ tend faiblement vers $\mu$.
\end{theoreme}
Le cas de dimension sup\'erieure nous semble beaucoup plus 
compliqu\'e. 
Lorsque $(F_n)$ est la suite des it\'er\'es d'une CM $F$ sur une vari\'et\'e
de grande dimension, pour d\'emontrer un r\'esultat analogue, on
fait appel \`a une propri\'et\'e de m\'elange de $F$ \cite{BriendDuval1,Dinh}. 
Dans le th\'eor\`eme 4.1,
l'hypoth\`ese (H4) est n\'ecessaire comme le montre 
l'example suivant.
\begin{exemple} \rm
Notons $z$ une coordonn\'ee affine de $\P^1$. Soient 
$F_n:\P^1\longrightarrow\P^1$, $F_n(z):=z^n+z$. La suite $(F_n)$ v\'erifie
(H1), (H2), (H3). La mesure $\mu$ est \'egale \`a la mesure 
de probabilit\'e invariante sur le cercle unit\'e. 
La mesure $S_\infty-\mu$ est \'egale \`a la masse de Dirac en $z=\infty$.
L'hypoth\`ese (H4) n'est
pas satisfaite. On v\'erifie que ces applications n'admettent aucun point fixe 
r\'egulier r\'epulsif.
\end{exemple}
\par
Le lemme suivant donne une majoration de nombre de points fixes au cas de 
dimension 1. Pour le cas de dimension sup\'erieure, on peut utiliser les arguments 
de la proposition 5.7. 
\begin{lemme} Soit $F$ une correspondance de graphe 
$\Gamma$ sur une surface de Riemann compacte $X$. Soient $d$
le degr\'e topologique et $\lambda$ le degr\'e interm\'ediaire d'ordre $0$ 
de $F$. Supposons que la diagonale $\Delta$ de $X\times X$ ne soit pas 
composante de $\Gamma$. Alors $F$ admet au plus $\lambda+d+2g\sqrt{\lambda d}$
points fixes, compt\'es avec multiplicit\'es, o\`u $g$ est le genre de $X$.
\end{lemme}
\begin{preuve} Observons que $d$ et $\lambda$ sont les degr\'es 
de $\Pi_{2|\Gamma}$ et $\Pi_{1|\Gamma}$. L'op\'erateur $F^*$ agit sur les 
groupes de Dolbeault $\H^{p,q}$ de $X$ (on identifiera $\H^{1,0}$ avec 
l'espace des $(1,0)$-formes holomorphes sur $X$). 
D'apr\`es la formule de Lefschetz 
\cite[p.314]{Fulton}, le nombre de points fixes de $F$ est \'egal \`a
\begin{eqnarray*}
p & = & \trace F^*_{|\H^{0,0}}+\trace F^*_{|\H^{1,1}} -
2\Re(\trace F^*_{|\H^{1,0}}) \\
& = & \lambda + d - 2\Re(\trace F^*_{|\H^{1,0}}).
\end{eqnarray*}
Soit $\alpha$ une $(1,0)$-forme 
holomorphe non-nulle qui est un vecteur propre de  $F^*_{|\H^{1,0}}$. 
Soit $\eta$ la valeur propre associ\'ee \`a $\alpha$.
Puisque $\deg \Pi_{1|\Gamma}=\lambda$, on a
\begin{eqnarray*}
F^*(\alpha\wedge\overline\alpha) & = & (\Pi_1)_* 
\big [(\Pi_{2|\Gamma})^* (\alpha)
\wedge (\Pi_{2|\Gamma})^* (\overline \alpha) \big]\\
& \geq &  \lambda^{-1}
\big[(\Pi_1)_* (\Pi_{2|\Gamma})^* (\alpha)\big] \wedge 
\big[(\Pi_1)_* (\Pi_{2|\Gamma})^* (\overline \alpha) \big]\\
& = & \lambda^{-1} F^*(\alpha)\wedge F^*(\overline\alpha) = \lambda^{-1}
|\eta|^2 (\alpha\wedge \overline\alpha).
\end{eqnarray*}
Or la masse de la mesure $F^*(\alpha\wedge\overline\alpha)$ est \'egale \`a 
$d$ fois celle de $\alpha\wedge \overline\alpha$. On en d\'eduit que 
$|\eta|\leq \sqrt{\lambda d}$. Donc 
$-\Re (\trace F^*_{|\H^{1,0}}) \leq g \sqrt{\lambda d}$ car $\dim \H^{1,0}=g$.
\end{preuve}
{\bf D\'emonstration du th\'eor\`eme 4.1.} 
L'hypoth\`ese (H1) implique que la multiplicit\'e de la diagonale $\Delta$ 
de $X\times X$ dans le graphe $\Gamma_n$  de $F_n$ 
est de l'ordre $\o(d_n)$. On peut 
donc supprimer $\Delta$ dans les $\Gamma_n$ pour simplifier la preuve.
\par
D'apr\`es le lemme 4.3, l'hypoth\`ese (H1)
implique que 
$\#\FR_n\leq d_n+\o(d_n)$. Par cons\'equent, quitte \`a extraire une sous-suite,
on peut supposer que $\widehat\mu_n$ converge vers une mesure $\widehat\mu$
de masse au plus \'egale \`a 1. 
Fixons un $\epsilon$, $0<\epsilon<1$. Il suffit de montrer que la masse de
$\mu-\widehat\mu$ est plus petite que $\epsilon$.
\par
Soit $\nu>0$ assez petit. Puisque $S_\infty$ est singuli\`ere par rapport \`a 
$\mu$, on peut trouver un ouvert $U$ tel que $S_\infty(U)<\nu$, 
$S_\infty(\partial U)=0$, $\mu(U)>1-\epsilon/4$ et $\mu(\partial U)=0$. 
En modifiant un peu l'ouvert $U$, on peut suposer qu'il est connexe et 
simplement connexe. Il est donc biholomorphe au disque unit\'e de $\C$.
Choisissons aussi un ouvert $U'\subset\subset U$ tel que
$\mu(U')>1-\epsilon/4$. 
\par
Quitte \`a diminuer l\'eg\`erement la taille de $U$ (lemme 3.6), on peut 
construire, comme au paragraphe 3, pour $n$ assez grand, 
$(1-\epsilon/4) d_n$ $F_n$-branches inverses
de $U$, qui sont de taille $\O(\lambda_n^{1/2}d_n^{-1/2})$. Notons 
$\B_{n,j}$ ces branches inverses et $F^{\B_{n,j}}$ les applications inverses 
associ\'ees.
\par
Fixons un point $x\in U\setminus\E$. On a $\mu^x_n\rightarrow \mu$. Notons
$\F_n$ la famille des $\B_{n,j}$ telles que $F^{\B_{n,j}}(x)\in U'$. 
Puisque $\mu(U')>1-\epsilon/4$, 
on a $\#\F_n> (1-\epsilon/2)d_n$ pour $n$ assez grand. 
De plus, comme $\B_{n,j}$ est de taille $\O(\lambda_n^{1/2}d_n^{-1/2})$, on a 
$F^{\B_{n,j}}(U)\subset\subset U$ pour $n$ grand. 
La derni\`ere relation implique que 
$F^{\B_{n,j}}$ admet un point fixe unique $x_{n,j}$ et sa d\'eriv\'ee
en $x_{n,j}$ est strictement inf\'erieure \`a 1 en module. C'est donc un point 
fixe r\'epulsif de $F_n$. 
\par
Posons $x_{n,j}':=F^{\B_{n,j}}(x)$. Quitte \`a extraire une sous-suite, 
on peut supposer que $d_n^{-1}\sum \delta_{x_{n,j}'}$ converge vers une mesure 
$\mu'$. Puisque $\#\F_n> (1-\epsilon/2)d_n$ et $\mu_n^x\rightarrow \mu$, 
la masse de la mesure positive $\mu-\mu'$ est inf\'erieure \`a $\epsilon/2$.
D'autre part, comme la taille de $\B_{n,j}$ tend vers 0, 
$d_n^{-1}\sum \delta_{x_{n,j}}$ tend aussi vers $\mu'$. Donc 
$\mu'\leq \widehat \mu$. On en d\'eduit que la masse de $\mu-\widehat\mu$ 
est major\'ee par $\epsilon$. La preuve du th\'eor\`eme 4.1 est achev\'ee.
\par
\hfill $\square$
\\
\section{Courants laminaires}
Soit $X$ une vari\'et\'e projective de dimension $k$.
Soit $(Z_n)$ une suite de sous-ensembles analytiques de dimension $s$ de $X$.
Nous cherchons des conditions suffisantes
pour que les valeurs adh\'erentes \`a la suite des courants 
d'int\'egration sur $Z_n$, proprement normalis\'es, soient 
(faiblement) laminaires.
Nous appliquons ce r\'esultat pour montrer que les courants de Green, de 
certaines
bidimensions, d'un automorphisme polynomial r\'egulier,
sont (faiblement) laminaires.
Rappelons d'abord la notion de laminarit\'e introduite par Bedford, Lyubich, 
Smillie  \cite{BedfordLyubichSmillie, Dujardin,deThelin}.
\par
Un courant positif ferm\'e $T$ de bidimension $(s,s)$ sur 
le polydisque $\Delta_1^k$ est appel\'e {\it lamin\'e} 
s'il s'\'ecrit $T=\int [\Gamma_a] \d \Lambda(a)$ o\`u  
les $\Gamma_a$ sont des graphes, deux \`a deux disjoints, 
d'applications holomorphes, 
au-dessus de $s$ directions du polydisque et o\`u $\Lambda$
est une mesure positive port\'ee par une transversale aux graphes. 
Un courant positif ferm\'e $T$ de bidimension $(s,s)$ sur 
une vari\'et\'e $V$ est appel\'e {\it uniform\'ement laminaire} si tout 
point 
$x\in V$ admet un voisinage biholomorphe \`a $\Delta_1^k$ sur lequel $T$ est 
lamin\'e.
Le courant $T$ est dit {\it laminaire} sur $V$ s'il existe une 
suite croissante de courants positifs $(T_i)$, des ouverts
$V_i\subset V$ de mesure $T_i\wedge \omega^s$ 
totale, tels que $T_i$ soit uniform\'ement 
lamin\'es dans $V_i$ pour tout $i$ et tels que $\lim T_i=T$.
Si, dans les d\'efinitions pr\'ec\'edentes, 
on ne suppose pas que les graphes sont 
disjoints, on dira qu'alors $T$ est {\it faiblement lamin\'e, uniform\'ement 
faiblement laminaire ou faiblement laminaire}.
\par
Soient $Y_n$ des vari\'et\'es projectives de dimension $s$, $1\leq s\leq k-1$.
Soient $\varphi_n:Y_n\longrightarrow X$ des applications holomorphes
\`a fibres discr\`etes et $Z_n$ leurs images.
Notons $m_n(z)$ le nombre de composantes irr\'eductibles locales de $Z_n$
en $z$. 
Nous faisons aussi l'hypoth\`ese que 
$m_n(z)=1$ pour $z\in Z_n$ g\'en\'erique.
On a $m_n(z):=\#\varphi_n^{-1}(z)$. 
Lorsque $s=k-1$, l'ensemble analytique 
$(m_n\geq 2)$ est de dimension pure $k-1$.
Dans ce cas, notons $\Sigma_n$ la $(k-1)$-cha\^{\i}ne holomorphe 
port\'ee par $(m_n\geq 2)$ dont les multiplicit\'es sont donn\'ees
par la fonction $m_n$. 
Notons finalement $v_n$ (resp. $c_n$) le volume de $Z_n$ (resp. de $\Sigma_n$) et $R$ 
l'ensemble des points adh\'erents \`a la suite des ensembles $(m_n\geq 2)$. 
\begin{theoreme} Supposons que $Y_n=\P^s$ pour tout $n$. Soit $T$ une valeur
adh\'erente de la suite de courants $(v_n^{-1}[Z_n])$. Alors $T$
est faiblement laminaire sur $X$. 
Si $s=k-1$, $T$ est laminaire sur
$X\setminus R$. 
Dans ce cas,  si la suite $(c_nv_n^{-1})$ 
est born\'ee ou si $T$ ne charge pas $R$, il   
est laminaire sur $X$.
\end{theoreme}
\par
Montrons d'abord la premi\`ere partie du th\'eor\`eme. 
Quitte \`a extraire une sous-suite, on peut supposer que 
la suite $(v_n^{-1}[Z_n])$ converge vers $T$. 
Si $(v_n)$ est born\'ee, $T$ est port\'e par un ensemble analytique 
de dimension $s$. On peut donc supposer que $(v_n)$ tend vers l'infini. 
Puisque $X$ est projective, on peut 
plonger $X$ dans un espace projectif. Donc
on peut supposer que $X=\P^k$. Fixons un sous-espace 
projectif $I$ de dimension $k-s-1$ de $\P^k$  tel
que $I\cap Z_n=\emptyset$ pour tout $n$.
Notons $\pi:\P^k\setminus I\longrightarrow \P^s$ la projection de centre $I$. 
On d\'efinit cette projection de la mani\`ere suivante. Fixons un sous-espace
projectif $\P^s$ dans $\P^k\setminus I$. Pour tout $z\in\P^k\setminus I$,
$\pi(z)$ est le point d'intersection de $\P^s$ avec le sous-espace projectif
contenant $I$ et passant par $z$.
On choisit $I$ de sorte que
$T$ ne charge pas $\overline{\pi^{-1}(Y)}$ pour
tout sous-ensemble analytique propre $Y$ de $\P^s$. 
Posons $\psi_n:=\pi\circ\varphi_n$. Ce sont des 
applications holomorphes de $\P^s$ dans $\P^s$.
Notons $d_n$ le degr\'e topologique de $\psi_n$. 
Observons que le degr\'e de $Z_n$ est aussi \'egal \`a $d_n$.
Puisque $v_n\sim d_n$, quitte \`a extraire une sous-suite de $(Z_n)$, 
on peut supposer  que $\lim v_n/d_n=c$.
\par
Soit $F_n:\P^k\longrightarrow \P^s$ la STM dont le graphe est l'ensemble
des points $(z,\pi(z))$ avec $z\in Z_n$. Le degr\'e topologique de $F_n$ 
est \'egal \`a $d_n$. Le degr\'e interm\'ediaire d'ordre $s-1$ est \'egal 
au degr\'e de $F_n^{-1}(D)$ pour une droite projective g\'en\'erique $D$ 
de $\P^s$. Il
est donc aussi \'egal \`a $d_n$. Le courant postcritique $S_n$ 
de $F_n$ est \'egal \`a celui de $\psi_n$. Par cons\'equent, sa masse est 
born\'ee par $s+1$ (\voir \cite{Sibony2} ou proposition 5.3). 
On peut donc supposer que la suite $S_n$ converge 
vers un courant $S_\infty$. On peut \'egalement supposer que la suite 
des courants $\Omega_n:=d_n^{-1}(F_n)_*(\omega_\FS)$
converge vers un courant 
$\Omega_\infty$. Fixons une constante $\eta>0$ assez petite.
\begin{lemme} Soit $x\in \P^s$ tel que $\nu(S_\infty,x)=0$
et $\nu(\Omega_\infty,x)=0$. Alors pour tout 
$\epsilon>0$, il existe $r>0$ tel que, si $n$ est 
assez grand, $\Ball(x,r)$ admette 
au moins $(1-\epsilon)d_n$ $F_n$-branches inverses de taille $\eta$.
\end{lemme}
\begin{preuve} Fixons une constante $\nu>0$ assez petite. 
D'apr\`es la d\'efinition du nombre de Lelong \cite{Skoda}, 
si $r_0>0$ est assez petit, 
les masses de $S_\infty$ et de
$\Omega_\infty$ sur $\Ball(x,2r_0)$ sont plus petites que $r_0^{2s-2}\nu$.
Pour $n$ assez grand, les masses de $S_n$ et de
$\Omega_n$ sur $\Ball(x,r_0)$ sont plus petites que $r_0^{2s-2}\nu$.
Puisque $\nu$ est petit, 
le lemme se d\'eduit directement de la proposition 3.2.
\end{preuve}
\par
On choisit une suite croissante de compacts $(K_i)$ dont la r\'eunion est 
\'egale \`a $\P^s\setminus L$ o\`u  
$L:=\{\nu(S_\infty,x)>0\}\cup \{\nu(\Omega_\infty,x)>0\}$. 
D'apr\`es un th\'eor\`eme de Siu \cite{Siu}, $L$ est une r\'eunion finie
ou d\'enombrable d'ensembles analytiques.
Rappelons qu'on a choisi $I$ tel que $T$ 
ne charge pas $\overline{\pi^{-1}(L)}$. 
On recouvre les $K_i$ par des familles finies 
de boules $\Ball_{i,j}$ qui v\'erifient
le lemme 5.2 pour $\epsilon=1/i$. 
Choisissons des ouverts connexes
$\Ball_{i,j}'\subset\Ball_{i,j}$ \`a bord lisse par morceaux tels que 
\begin{enumerate}
\item Pour chaque $i$, les ouverts $\Ball_{i,j}'$ sont disjoints;
\item Chaque ouvert $\Ball_{i+1,j}'$ est contenu dans l'un des $\Ball'_{i,j}$;
\item $T$ ne charge pas $\pi^{-1}(\partial \Ball_{i,j}')$;
\item La r\'eunion $\Ball_i':=\cup_j\Ball_{i,j}'$ v\'erifie $\overline \Ball_i' 
\supset K_i$.
\end{enumerate}
Consid\'erons les 
$F_n$-branches inverses de $\Ball_{i,j}'$ qui sont de taille $\eta$ 
(\voir lemme 5.2). 
Notons $T_{n,i}$ la somme (en $j$) 
des courants d'int\'egration sur les images de ces $F_n$-branches inverses.
Observons que, puisque $\eta$ 
est petit, 
la famille des applications holomorphes de $\Ball_{i,j}'$ dans $\P^k$, 
dont les images sont de taille $\eta$, est normale.
Quitte \`a extraire une sous-suite, on peut supposer que 
$v_n^{-1}T_{n,i}$ converge, dans $\pi^{-1}(\Ball_i')$, vers 
un courant uniform\'ement faiblement laminaire $cT_i$  
quand $n\rightarrow\infty$.
Le choix des ouverts $\Ball_{i,j}'$ implique que la suite 
$(T_i)$ est croissante.
Comme elle est domin\'ee par $cT$, elle converge vers un courant 
faiblement laminaire $cT_\infty$. Par construction,
on a $T\geq T_\infty$ et $(T-T_\infty)\wedge\pi^*(\omega_\FS^s)=0$ sur 
$\P^k\setminus I$. 
\par
Soient $\pi'$ une autre projection g\'en\'erique de centre $I'$. 
On construit de la
m\^eme mani\`ere un courant faiblement laminaire 
$T_\infty'$. Observons que, par construction,  $T_\infty'$ domine
les $T_i$ et donc $T_\infty'\geq T_\infty$. 
Par sym\'etrie, $T_\infty=T'_\infty$. On en d\'eduit que 
$(T-T_\infty)\wedge(\pi')^*(\omega_\FS^s)=0$ sur 
$\P^k\setminus I'$ pour tout $\pi'$ g\'en\'erique. Donc $T$ est \'egal 
au courant faiblement laminaire $T_\infty$. 
\par
\
\par
Supposons maintenant que $s=k-1$. 
On peut consid\'erer
$X$ comme une sous-vari\'et\'e d'un espace projectif $\P^{k'}$. 
Soient $V$ et $V'$  des voisinages de $R$ avec $V'\subset\subset V$.
On construit les $F_n$-branches inverses et les courants $T_{n,i}$, $T_i$ 
comme ci-dessus pour un $\eta$ plus petit que
$\dist(\partial V,\partial V')$. 
\par 
Consid\'erons uniquement les branches inverses de $\Ball_{i,j}'$ 
dont les images intersectent $X\setminus \overline V$. 
Les images des $\Ball_{i,j}'$  sont des morceaux de vari\'et\'e complexe
de dimension $k-1$ 
de diam\`etre plus petit que $\dist(\partial V,\partial V')$.
Elles ne rencontrent pas $V'$. 
Notons $T_{n,i}'$ la somme 
des courants d'int\'egration sur ces morceaux de vari\'et\'e. 
On a $T_{n,i}'=T_{n,i}$ sur $X\setminus\overline V$. Par d\'efinition 
de $R$, lorsque $n$ est assez grand, 
les morceaux de vari\'et\'e d\'efinissant $T_{n,i}'$ sont disjoints. 
Dans $\pi^{-1}(\Ball_i')\setminus \overline V'$,  puisque $s=k-1$, 
les diff\'erentes limites de suites de tels morceaux 
sont aussi disjointes. On en d\'eduit que $T_i$ est uniform\'ement laminaire 
dans $\pi^{-1}(\Ball_i')\setminus \overline V$. 
Par suite, $T$ est laminaire dans $X\setminus (I\cup \overline V)$
et donc dans $X\setminus R$. Si $T$ ne charge pas $R$, il est laminaire 
dans $X$.
\par
Supposons maintenant que la suite $(c_nv_n^{-1})$ soit born\'ee.
On peut donc supposer que la suite des courants $d_n^{-1}[\Sigma_n]$ 
converge vers un courant $S'_\infty$.
On peut appliquer 
la proposition 3.8 pour 
des petites boules dont les centres n'appartiennent pas
\`a $\{\nu(\pi_*S'_\infty,z)>0\}$. Ceci nous permet de construire 
des $F_n$-branches inverses de $\Ball_{i,j}'$ qui ne rencontrent pas 
$|\Sigma_n|$ et donc qui ne s'intersectent pas. On obtient, 
par cons\'equent, que $T$ est
laminaire sur $X$. 
\par
\hfill $\square$
\par
Le th\'eor\`eme 5.1 est aussi valable pour les deux cas suivants o\`u on n'a 
pas n\'ecessairement $Y_n=\P^s$:
\begin{enumerate}
\item Le fibr\'e canonique $\KK_{Y_n}$ de $Y_n$ est n\'egatif pour tout $n$;
\item Le groupe de Picard $\Pic(Y_n)$ 
de $Y_n$ est isomorphe \`a $\Z$ et la classe de
Chern $c_n$ de $\KK_{Y_n}$ v\'erifient $\int c_n^s=\O(v_n)$. 
\end{enumerate}
Dans le cas de dimension 1, la condition 2 ci-dessus \'equivaut au fait que le genre 
de $Y_n$ est de l'ordre au plus \'egal \`a $\O(v_n)$ \cite[p.216]{GriffithsHarris}.
\par
Pour la d\'emonstration, il suffit de montrer que la masse du courant
postcritique associ\'e \`a l'application $\psi_n:=\pi\circ\varphi_n$ est 
de l'ordre $\O(v_n)$.
La proposition suivante donne l'estimation n\'ecessaire.
\begin{proposition} Soient $Y$ une vari\'et\'e projective de dimension $s$ 
et $c$ la classe de Chern de la fibr\'e canonique $\KK_Y$ de $Y$. 
Soient $\psi:Y\longrightarrow \P^s$ une application holomorphe \`a fibres 
discr\`etes, $d$ son degr\'e topologique et $S$ son courant postcritique. 
Alors, si $\KK_Y$ est n\'egatif, la masse de $S$ est major\'ee par $(s+1)d$.
Si $\Pic(Y)\simeq \Z$ et si $\KK_Y$ est positif, 
la masse de $S$ est \'egale \`a 
$(s+1)d+d^{1-1/s}\big(\int c^s\big)^{1/s}$. 
\end{proposition} 
\begin{preuve} 
Notons $R$ le diviseur de 
ramification de $\psi$. Le courant $S$ est \'egal au courant 
d'int\'egration sur 
$\psi_*(R)$. 
Dans la suite, la notation $[\ ]$ d\'esigne les classes de cohomologie. 
Soit $H$ un hyperplan de $\P^s$. La classe de Chern de $\KK_{\P^s}$ est 
\'egale \`a $-(s+1)[H]$ \cite[p.146]{GriffithsHarris}. 
D'apr\`es la formule de Riemann-Hurwitz \cite[p.62]{Fulton}, on a
$[R]=(s+1)\psi^*[H]+c$. Puisque $\psi_*\circ\psi^*=d.\id$, on a 
$[S]=(s+1)d[H] + \psi_*(c)$. On en d\'eduit que si $\KK_Y$ est n\'egatif,
$[S]\leq (s+1)d[H]$. Donc la masse de $S$ est major\'ee par $(s+1)d$.
\par
Supposons maintenant que $\Pic(Y)\simeq \Z$ et $\KK_Y$ est positif. 
Il existe une constante positive $\lambda$
telle que $c=\lambda\psi^*[H]$. 
On a $[S]=(s+1)d[H]+d\lambda [H]$. D'autre part,
puisque $\deg \psi=d$, on a $\int\psi^* [H]^s=d$. 
Donc $\int c^s=d\lambda^s$.
La proposition s'en d\'ecoule.
\end{preuve}
\begin{remarque} \rm
Dans le cas de dimension $s=1$, les r\'esultats ci-dessus sont 
d\'emontr\'es par Bedford-Lyubich-Smillie \cite{BedfordLyubichSmillie}, 
Dujardin \cite{Dujardin} et 
ils sont aussi valables quand $X$ n'est pas une vari\'et\'e projective, \voir 
de Th\'elin \cite{deThelin}. 
\end{remarque}
\par
Nous allons maintenant donner une hypoth\`ese plus faible sur la suite
$(Z_n)$ telle que le th\'eor\`eme 5.1 reste vrai. Ceci n\'ecessite 
l'introduction des objets assez sophistiqu\'es. Nous supposons pour
simplifier que les $Z_n$ sont lisses mais la m\'ethode permet 
aussi de traiter le cas 
des vari\'et\'es singuli\`eres.
\par
Soit $TX$ le fibr\'e tangent holomorphe de $X$. Soit $\Lambda^{k-s} (TX)$ 
le fibr\'e
des $(k-s,0)$-vecteurs tangents holomorphes. Notons $X_{k-s}$ 
la projectivisation
de  $\Lambda^{k-s}(TX)$ et $\Pi_{k-s}:X_{k-s}
\longrightarrow X$ la projection canonique.
Observons que les fibres de $\Pi_{k-s}$ sont des espaces 
projectifs de m\^eme dimension. 
Soit $\widehat X_{k-s}$ l'ensemble des points 
$[x,v]\in X_{k-s}$
avec $x\in X$ et $v$ un $(k-s,0)$-vecteur tangent simple \cite[pp.23-4]{Federer}.
C'est une sous-vari\'et\'e de $X_{k-s}$. 
Notons $\widehat Z_n$ l'ensemble des points
$[x,v]\in \Pi_{k-s}^{-1}(Z_n)
\cap \widehat X_{k-s}$ tels que le $(k-s,0)$-vecteur 
$v$ ne soit pas transversal \`a la vari\'et\'e $Z_n$. C'est une 
sous-vari\'et\'e de $\Pi_{k-s}^{-1}(Z_n)\cap \widehat X_{k-s}$. Observons que 
les fibres de $\Pi_{k-s|\widehat Z_n}$ (resp. $\Pi_{k-s|\widehat X_{k-s}}$) 
sont des vari\'et\'es (resp. des grassmanniennes) identiques. La fibre 
$\Pi_{k-s}^{-1}(x)\cap \widehat X_{k-s}$ param\`etre les sous-espaces complexes de 
dimension $k-s$ de l'espace tangent de $X$ en $x$. 
\begin{theoreme} Soient $Z_n$ des sous-vari\'et\'es complexes
de dimension $s$ de volume $v_n$ 
d'une vari\'et\'e projective $X$ de dimension $k$. 
Supposons que les volumes $\widehat v_n$ de $\widehat Z_n$ v\'erifient
$\widehat v_n =\O(v_n)$. Alors toute valeur adh\'erente \`a la suite 
$(v_n^{-1}[Z_n])$ est un courant faiblement laminaire sur $X$. 
C'est un courant laminaire si $s=k-1$.
\end{theoreme}
\begin{preuve} On reprend la d\'emonstration du th\'eor\`eme 5.1.
Puisque $X$ est projective, on peut la plonger dans un espace projectif.
Pour simplifier la preuve, 
on peut supposer $X=\P^k$. Notons $G$ la grassmannienne qui 
param\`etre les sous-espaces projectifs de dimension $k-s$ de $\P^k$. 
On associe tout $(k-s,0)$-vecteur simple $v$, tangent \`a $\P^k$ 
en $x$, au 
sous-espace projectif de dimension $k-s$ passant par $x$ et tangent \`a $v$.
Ceci d\'efinit une application holomorphe $\Phi$ de $\widehat X_{k-s}$ 
dans $G$.
Posons $\widetilde Z_n:=\Phi(\widehat Z_n)$. Cet ensemble param\`etre les 
sous-espaces projectifs de dimension $k-s$ qui n'intersectent 
pas $Z_n$ transversalement. L'hypoth\`ese sur $\widehat v_n$ implique que
les volumes $\widetilde v_n$ des $\widetilde Z_n$ v\'erifient $\widetilde 
v_n=\O(v_n)$. Dans les notations du th\'eor\`eme 5.1, ceci implique
que les courants $S_n$ sont de masse uniform\'ement born\'ee
pour $\pi$ g\'en\'erique. 
En effet, l'ensemble des fibres de $\pi$ correspond \`a une sous-vari\'et\'e
$V_\pi$ de $G$. C'est le degr\'e de $V_\pi\cap \widetilde Z_n$ 
qui donne l'estimation 
de masse pour $S_n$.  
Le th\'eor\`eme 5.5 se d\'emontre exactement comme le th\'eor\`eme 5.1.
\end{preuve}
\par
Consid\'erons \`a pr\'esent un automorphisme polynomial $f$ de $\C^k$.
On peut prolonger cet automorphisme en une application m\'eromorphe de $\P^k$ 
dans $\P^k$. Notons aussi $f$ et $f^{-1}$ les prolongements de $f$ et 
$f^{-1}$. Soient $d_\pm$ et $I_\pm$ le degr\'e alg\'ebrique et 
l'ensemble d'ind\'etermination 
de $f^{\pm 1}$. On dit que $f$ est {\it r\'egulier} (au sens de Sibony)
si $I_+\cap I_-=\emptyset$. Dans le cas de dimension $k=2$, ce sont les
automorphismes du type H\'enon. Sibony a montr\'e \cite{Sibony2} 
qu'il existe un $s$, $1\leq s\leq k-1$ tel que $\dim I_+=s-1$, 
$\dim I_-=k-s-1$ et $d_+^{k-s}=d_-^s$. Les ensembles d'ind\'etermination de 
$f^{\pm n}$ sont aussi \'egaux \`a $I_\pm$.
Sibony a construit les courants de Green 
$T_\pm$ de bidegr\'e $(1,1)$ associ\'es \`a $f^{\pm 1}$ et montr\'e que 
les produits ext\'erieurs 
$T_+^p\wedge T_-^q$, avec $0\leq p\leq k-s$ et $0\leq q\leq s$,
d\'efinissent des courants non-nuls, positifs, ferm\'es et invariants par $f$. 
Ces courants ne chargent pas les sous-ensembles analytiques propres de $\P^k$.
D'apr\`es un th\'eor\`eme de Russakovskii-Shiffman
\cite{RussakovskiiShiffman, DinhSibony3}, si $H_\pm$ sont des sous-espaces
projectifs g\'en\'eriques de dimensions respectives $s$ et $k-s$, les suites
de courants $d_+^{-(k-s)n}(f^n)^*[H_+]$ et $d_-^{-sn}(f^{-n})^*[H_-]$  
tendent faiblement 
vers $T_+^{k-s}$ et $T_-^s$. Le r\'esultat suivant g\'en\'eralise
un th\'eor\`eme de Bedford-Lyubich-Smillie \cite{BedfordLyubichSmillie}.
\begin{corollaire} Les courants $T_+^{k-s}$ et $T_-^s$ sont 
faiblement laminaires. Si $s=k-1$, le courant $T_+$ est laminaire,
si $s=1$, le courant $T_-$ est laminaire. 
\end{corollaire}
\begin{preuve} Soit $H_+$ un sous-espace g\'en\'erique 
de dimension $s$ 
qui n'intersecte pas $I_-$. 
Les applications $f^{-n}$ sont holomorphes sur $H_+$.
Consid\'erons la suite d'ensembles analytiques $Z_n:=f^{-n}(H_+)$. 
On a $T_+^{k-s}=\lim d_+^{-(k-s)n}[Z_n]$. 
Le th\'eor\`eme
5.1 implique que $T_+^{k-s}$ est faiblement laminaire.
Si $s=k-1$, l'ensemble $R$ associ\'e \`a la suite $(Z_n)$ est contenu dans 
l'hyperplan \`a l'infini qui n'est pas charg\'e par $T_+$. 
Donc $T_+$ est laminaire. La preuve est identique pour $T_-^s$.
\end{preuve}
\par
Consid\'erons maintenant un automorphisme holomorphe $f$ d'une vari\'et\'e 
projective complexe $X$ de dimension $k$. Pour tout $0\leq p\leq k$, 
il existe une constante $d_p\geq 1$ et un entier $l_p\geq 0$ tels que 
la norme de $f^{n*}$ sur le groupe de Dolbeault $\H^{p,p}(X,\C)$ 
soit de l'ordre $n^{l_p}d_p^n$. Ceci se voit ais\'ement avec la forme de Jordan
de la matrice associ\'ee \`a $f^*$. On dira que $d_p$ est 
{\it le degr\'e dynamique d'ordre $p$} de $f$. On a $d_0=d_k=1$ et $l_0=l_k=0$. 
D'apr\`es un th\'eor\`eme 
de Khovanskii-Teissier, la suite $(d_p)$ est concave \cite{Guedj}. 
En particulier, si $d_p>d_{p+1}$, on a $d_p>\max_{q\geq p+1} d_q$.
\par
Observons que si $Z$ est
un sous-ensemble analytique de dimension $k-p$ de $X$, le volume de
$f^{-n}(Z)$ est de l'ordre au plus \'egal \`a $n^{l_p}d_p^n$. Pour 
estimer les normes de $f^{n*}$ sur $\H^{p,q}(X,\C)$, nous avons la proposition 
suivante.
\begin{proposition} Il existe $c>0$ tel que pour tout $n\geq 1$ la norme 
$A_{p,q,n}$ de $f^{n*}$ sur $\H^{p,q}(X,\C)$
v\'erifie 
$$A_{p,q,n}^2\leq cn^{l_p+l_q}d_p^nd_q^n.$$
En particulier, le rayon spectral $r_{p,q}$ 
de $f^*$ sur $\H^{p,q}(X,\C)$ v\'erifie 
$$r_{p,q}\leq \sqrt{d_pd_q}.$$
\end{proposition}
\begin{preuve} 
Fixons une forme de K\"ahler $\omega$ sur $X$.
Rappelons que d'apr\`es la th\'eorie de Hodge \cite[p.116]{GriffithsHarris},  
$\H^{p,q}(X,\C)\simeq
\overline{\H^{q,p}(X,\C)}$.
On en d\'eduit que $A_{p,q,n}= A_{q,p,n}$.
On peut donc supposer $p\leq q$.
Soit $\varphi$ (resp. $\psi$) une $(p,q)$-forme 
(resp. $(q,p)$-forme) lisse $\overline\partial$-ferm\'ee sur $X$. 
Consid\'erons  
l'automorphisme $F$ de $X\times X$ d\'efini par $F(x,y)=(f(x),f(y))$.
Posons $\Phi:=\varphi(x)\wedge \psi(y)$ et $\Omega:=\omega(x)+\omega(y)$.
Il suffit de majorer la norme de la classe 
$F^{n*}[\Phi]$ dans $\H^{p+q,p+q}(X\times X,\C)$.
Fixons une $(2k-p-q,2k-p-q)$-forme positive ferm\'ee $\Theta$ sur $X\times X$.
Par la dualit\'e  de Poincar\'e \cite[p.53]{GriffithsHarris}, 
il faut seulement estimer l'int\'egrale $\int F^{n*}(\Phi)\wedge \Theta$.
Observons que $\Phi$ s'\'ecrit localement comme combinaison $\C$-lin\'eaire
de $(p+q,p+q)$-formes positives qui sont faiblement major\'ees par 
$\omega^p(x)\wedge \omega^p(y)\wedge (\omega^{q-p}(x)+\omega^{q-p}(y))$. 
D'autre part, $\Theta$ est aussi major\'ee par un multiple de 
$\Omega^{2k-p-q}$. D'o\`u
\begin{eqnarray*}
\left|\int_{X\times X} F^{n*}(\Phi)\wedge \Theta\right| 
& \lesssim &  \int_{X\times X}  
f^{n*}\omega^p(x)\wedge f^{n*}\omega^p(y)\wedge \\
& & \hspace{1cm}\wedge
\left(f^{n*}\omega^{q-p}(x)+f^{n*}\omega^{q-p}(y)\right)
\wedge \Omega^{2k-p-q} \\
& \lesssim & \left(\int_X f^{n*}\omega^p\wedge 
\omega^{k-p}\right)
\left(\int_X f^{n*}\omega^q\wedge \omega^{k-q}\right)\\ 
& \lesssim & n^{l_p+l_q}d_p^nd_q^n.
\end{eqnarray*}
La proposition en d\'ecoule.
\par
Nous rappelons ici la notion de positivit\'e faible utilis\'ee ci-dessus. 
Une $(s,s)$-forme $\Psi$ sur $X$ est 
{\it positive} si, pour tout $z\in X$, $\Psi(z)$ s'\'ecrit comme somme
de $(s,s)$-vecteurs de la forme
$$i\alpha_1 \wedge 
\overline \alpha_1\wedge \ldots \wedge i\alpha_s \wedge \overline \alpha_s$$
o\`u les $\alpha_j$ sont des 
$(1,0)$-vecteurs cotangents de $X$ en $z$.
On dira que $\Psi$ est 
{\it faiblement positive} 
si en tout point $z\in X$ et pour tous $(1,0)$-vecteurs cotangents 
$\alpha_j$ on a 
$$\Psi\wedge i\alpha_1\wedge \overline \alpha_1 \wedge \ldots \wedge i\alpha_{k-s}
\wedge \overline \alpha_{k-s}\geq 0.$$
Toute $(s,s)$-forme lisse sur $X$ s'\'ecrit comme combinaison lin\'eaire,
\`a coefficients dans ${\cal C}^\infty(X)$, de formes positives.
On dit que $\Psi$ est {\it faiblement domin\'ee} par $\Psi'$ 
si $\Psi'\pm \Psi$ sont faiblement positives.
Notons $(z_1,\ldots, z_k)$ des coordonn\'ees locales de $X$ en $z$ et 
$\d z_I:=\d z_{i_1}\wedge\ldots \wedge \d z_{i_s}$ pour tout multi-indice
$I=(i_1,\ldots,i_s)$. 
On peut v\'erifier que 
les parties r\'eelle et imaginaire de $\d z_I \wedge \d \overline z_J$ 
sont faiblement domin\'ees par des combinaisons de 
$i^{s^2}\d z_I\wedge \d \overline z_I$ 
et de $i^{s^2}\d z_J\wedge \d \overline z_J$.
\par
La forme $\Phi$ sur $X\times X$, consid\'er\'ee pr\'ec\'edemment, s'\'ecrit localement 
comme combinaison lin\'eaire, \`a coefficients dans ${\cal C}^\infty(X\times X)$, 
de formes du type 
$\d x_I \wedge \d \overline x_{I'} \wedge \d y_{J'} \wedge \d \overline y_J$ avec 
$|I|=|J|=p$ et $|J|=|J'|=q$. Ses parties r\'eelle est imaginaire sont donc 
faiblement major\'ees par un multiple de 
$\omega^p(x)\wedge \omega^p(y)\wedge (\omega^{q-p}(x)+\omega^{q-p}(y))$. 
\end{preuve}
\begin{remarque} \rm
Lorsque l'ensemble des points p\'eriodiques de $f$ n'a pas de composante de 
dimension strictement positive, la proposition ci-dessus permet d'estimer
le nombre de points p\'eriodiques d'ordre $n$ en fonction de $n$. Il suffit 
d'appliquer la formule de Lefschetz 
\cite[p.314]{Fulton} et ceci est valable aussi pour le cas o\`u $f$ est une 
application m\'eromorphe.
\end{remarque}
\begin{corollaire} Soit $f$ un automorphisme holomorphe sur une vari\'et\'e
projective $X$ de dimension $k$ comme ci-dessus. Soit $Z$ 
une sous-vari\'et\'e de dimension $s$ de $X$. Supposons que 
$d_{k-s}>d_{k-s+1}$. Alors toute valeur adh\'erente de la suite 
$(n^{-l_{k-s}} d_{k-s}^{-n}[f^{-n}(Z)])$ est un courant faiblement laminaire.
C'est un courant laminaire si $s=k-1$.
\end{corollaire}
\begin{preuve} On peut supposer que le volume de $f^{-n}(Z)$ est de l'ordre 
$n^{l_{k-s}}d_{k-s}^n$. Autrement, le courant limite est nul.
Posons $Z_n:=f^{-n}(Z)$. D'apr\`es le th\'eor\`eme 5.5, il suffit de 
montrer que le volume de $\widehat Z_n$ est de l'ordre 
$n^{l_{k-s}}d_{k-s}^n$.
\par
Rappelons que les fibres de $\Pi_{k-s}$ sont isomorphes \`a un espace projectif
$\P^N$. Notons $r$ la dimension de $\widehat Z_n$. 
On a $r=(k-s+1)s\geq k$ et $\dim X_{k-s}=N+k$.
L'application $f$ induit un automorphisme $F$ sur $X_{k-s}$ qui 
pr\'eserve la fibration $\F$ donn\'ee par $\Pi_{k-s}$. 
Soit $\Omega$ une $(r,r)$-forme strictement positive et ferm\'ee sur 
$X_{k-s}$. Il faut estimer $\int [\widehat Z_n]\wedge \Omega$. 
Posons $\Omega_n:=(F^n)_*\Omega$.
On a 
 $$\int [\widehat Z_n]\wedge \Omega =\int (F^n) ^*[\widehat Z]\wedge \Omega
=\int \Omega_n\wedge [\widehat Z].$$
\par
Observons que $\H^{p,q}(\P^N,\C)$ est nul si $p\not=q$ et 
$\H^{p,p}(\P^N,\C)\simeq \C$. Le groupe $\H^{2p}(\P^N,\Z)$ est isomorphe 
\`a $\Z$ et est engendr\'e par la classe d'un sous-espace projectif 
de dimension $N-p$.
Ceci implique que l'action du groupe $\pi_1(X)$ sur les  
cohomologies des fibres de $\F$, est triviale.
L'action de $F$, elle aussi, est triviale sur les  
cohomologies des fibres de $\F$.
D'apr\`es la th\'eorie de suites spectrales de Leray et un th\'eor\`eme de
Deligne
\cite[pp.438-68]{GriffithsHarris}, il existe un morphisme bijectif $\tau$
qui commute avec les actions de $f$ et de $F$
$$\tau:\ \H^{2r}(X_{k-s},\C)\longrightarrow 
\sum_{0\leq p,q\leq k\atop p+q\ {\rm pair}} \H^{p,q}(X,\C).$$
Le membre a droite est isomorphe \`a
$$\sum_{0\leq p,q\leq k\atop p+q\ {\rm pair}}
\H^{p,q}(X,\C) \otimes_\C \H^{2r-p-q}(\P^N,\C).$$
Notons $\tau_{p,q}$ la projection de $\H^{2r}(X_{k-s},\C)$ sur 
$\H^{p,q}(X,\C)$. Posons aussi $c_{p,q}:=\tau_{p,q}[\Omega]$
et $c_{p,q,n}:=\tau_{p,q}[\Omega_n]$. On a 
$c_{p,q,n}=(f^n)_*(c_{p,q})$. Il reste \`a estimer
$\int \tau^{-1}(c_{p,q,n}) \wedge [\widehat Z]$. 
\par
La fibration $\F$ est localement isomorphe \`a un produit. Si $U$ est un ouvert assez
petit de $X$, on a $\Pi_{k-s}^{-1}(U)\simeq U\times \P^N$. La classe 
$\tau^{-1}(c_{p,q,n})$ peut-\^etre repr\'esent\'ee par une forme qui s'\'ecrit 
localement comme somme de formes
du type $\eta\wedge \psi$ o\`u $\eta$ est une forme sur  $\Pi_{k-s}^{-1}(U)$ et 
o\`u $\psi$ est une forme sur $U$ de bidegr\'e $(p,q)$. Par cons\'equent,
lorsque $\max(p,q)>s$, on a  
$\int \tau^{-1}(c_{p,q,n}) \wedge [\widehat Z]=0$ car $\Pi_{k-s}(\widehat Z)$
est de dimension $s$.
\par
Il suffit maintenant de majorer $(f^n)_*c_{p,q}$ 
pour $p\leq s$ et $q\leq s$. Observons que la norme de $f_*$ sur 
$\H^{p,q}(X,\C)$ est \'egale \`a celle de $f^*$ sur $\H^{k-p,k-q}(X,\C)$.
D'apr\`es la proposition 5.7, on a 
$\|(f^n)_*c_{p,q}\|\lesssim n^{l_{k-s}}d_{k-s}^n$ car $d_{k-s}>d_m$ pour tout 
$m>k-s$.
La preuve du corollaire est achev\'ee.
\end{preuve}
\par
Le corollaire 5.9 n'affirme pas que la suite 
$(n^{-l_{k-s}} d_{k-s}^{-n}[f^{-n}(Z)])$ converge. Nous renvoyons le lecteur 
\`a \cite{DinhSibony3,DinhSibony4} pour 
l'\'etude de la convergence de cette suite  
vers des courants invariants ainsi
que pour les propri\'et\'es dynamiques de courants limites (courants de Green).
\par
\
\\
{\bf Remerciement.} 
Une partie de cet article a \'et\'e \'ecrite, avec le support 
d'Alexander von Humboldt Stiftung, pendant la visite de l'auteur   
\`a Humboldt-Universit\"at zu Berlin. Il tient \`a remercier
ces organisations ainsi que Professeur J\"urgen Leiterer pour
leur aide et leur accueil.
Tien-Cuong Dinh,
Math\'ematique - B\^at. 425, UMR 8628, 
Universit\'e Paris-Sud, 91405 Orsay, France. 
E-mails: Tiencuong.Dinh@math.u-psud.fr.
\end{document}